\documentstyle[12pt,bezier]{article}

\title{Global strong solution to Maxwell-Dirac equations in $1+1$ dimensions}
\author{Aiguo You and Yongqian Zhang  \\ {\small School of Mathematical Sciences,
Fudan University,
Shanghai 200433} \\ {\small Key Laboratory of Mathematics for
Nonlinear Sciences }
\\ {\small Email address: yongqianz@fudan.edu.cn} }
\begin{document}
\maketitle

\makeatletter
\renewcommand{\theequation}{\thesection.\arabic{equation}}
\@addtoreset{equation}{section}
\makeatother

\newtheorem{lemma}{Lemma}[section]
\newtheorem{proposition}{Proposition}[section]
\newtheorem{theorem}{Theorem}[section]
\newtheorem{definition}{Definition}[section]
\newtheorem{remark}{Remark}[section]
\def\R{\mathbb{R}}

\begin{abstract}
  The Maxwell-Dirac equations  with nonzero charge mass in one space dimension are studied under the Lorentz gauge condition. The global existence and uniqueness of solution in $C([0,+\infty);L^2(R^1))\times C(R^1\times [0,\infty))$ for initial value problem of Maxwell-Dirac equations are proved.
\end{abstract}

\section{Introduction}

We are concerned with the Maxwell-Dirac equations in one space dimension
\begin{equation}\label{eq:1}
\left\{\begin{array}{l}
(i {\gamma}^{\mu} D_{\mu}+ m I) {\Psi} = 0, \\
{\partial}_{\mu} F^{\mu\nu} = J ^{\nu}, \end{array}
\right.
\end{equation}
under the Lorentz gauge condition \[ \partial_{t} A_0 - \partial_{x} A_1 = 0, \]
with
\begin{equation}\label{eq:1-initialdata}
\left\{ \begin{array}{l}
\Psi_{j}(x,0) = \psi_{j}(x) \,  \, \, (j=1,2),  \\
A_{\nu}(x,0) = a^{0}_{\nu}(x)\,  \, \, (\nu=0,1), \\
\partial_{t} A_{\nu}(x,0) = a^{1}_{\nu}(x)\,  \, \, (\nu=0,1),
\end{array}
\right.
\end{equation}
which satisfy the constraint $a^1_{0}-\partial_{x} a^0_1 = 0.$
Here $ D_{\mu} = {\partial}_{\mu} -i A_{\mu} $ is the covariant derivative
and $F^{\mu\nu} = {\partial}_{\nu} A_{\mu} - {\partial}_{\mu} A_{\nu}$ is the curvature
associated with the gauge field $A_{\mu} \in R$. $ {\Psi} $ denotes a two-spinor field defined on $ R^{1+1} $, $J^\nu = \Psi ^{\dag} \gamma^{0} \gamma^{\nu}\Psi$ is a current density and $\Psi^{\dag} = (\overline{\Psi_1},\overline{\Psi_2})$
denotes the complex conjugate transpose of $\Psi$. ${\partial}_0 = {\partial}_t$,
${\partial}_1 = {\partial}_x$.

  The Dirac gamma matrices are of the following:
\begin{displaymath}
\gamma^0 =
\left( \begin{array}{ccc}
1 & 0  \\
0 & -1
\end{array} \right),
\gamma^1 =
\left( \begin{array}{ccc}
0 & 1  \\
-1 & 0
\end{array} \right).
\end{displaymath}

 The system (\ref{eq:1}) can be rewritten as follows,
\begin{equation} \label{eq:2}
\left\{\begin{array}{l}
\partial_{t} \Psi_1 + \partial_{x} \Psi_2 = im \Psi_1 + i A_{0}\Psi_1 + i A_{1}\Psi_2, \\
\partial_{t} \Psi_2 + \partial_{x} \Psi_1 = -im \Psi_2 + i A_{0}\Psi_2 + i A_{1}\Psi_1,\\
\Box A_{0} = |\Psi_1|^2 + |\Psi_2|^2,\\
\Box A_{1} = -(\overline{\Psi_{2}} \Psi_1 + \overline{\Psi_{1}} \Psi_2).
\end{array}\right.
\end{equation}

The global well-posedness of classical solution for the Maxwell-Dirac system in
$ R^{1+1} $ has been established in \cite{C}.
Many works are devoted to study the solution in different functional spaces since then, see for instance \cite{B}-\cite{Okamoto} and the references therein.
  Recently in \cite{H} Huh proved the global
 well-posedness of the strong solutions for the Maxwell-Dirac system(\ref{eq:1}) in
 $ R^{1+1} $, where he assumed that the mass of charge is zero, that is,
 $m=0$ in (\ref{eq:1}), and the solutions can be obtained by the explicit formula. For $m>0$, as far as we know, there is no explicit formula for solutions.
In this paper we consider more general case than that in \cite{H},  that is, $m\ge 0$.
Our aim is to find the $(\Psi, A)$ which solves the equation (\ref{eq:2}) with initial
data (\ref{eq:1-initialdata}) in the following sense and prove its uniqueness.
\begin{definition}\label{def-strongsolu} $(\Psi, A) $ with
$\Psi\in L^2_{loc}(R^1\times (0,+\infty))$ and $A\in C(R^1\times [0,+\infty))$ is called a strong solution to (\ref{eq:2})
with  the initial data (\ref{eq:1-initialdata}) provided that there exists a sequence of
smooth solutions $\{ (\Psi^n, A_{\mu}^n)\}_{n=1}^{\infty}$ such that
\[ ||\Psi^n-\Psi||_{L^2(K)}+||A^n_1-A_1||_{C(K)}+||A^n_2-A_2||_{C(K)} \to 0 \]
as $n \to \infty$ for any compact set $K\subset R^1\times [0,\infty)$,
and
\[ ||\Psi^n(\cdot,0)-\psi||_{L^2([a,b])}+||A^n_{\mu}(\cdot,0)-a_{\mu}^0||_{C([a,b])}+
||\partial_t A^n_{\mu}(\cdot,0)-a_{\mu}^1||_{C([a,b])} \to 0 \]
as $n\to \infty$ for $\mu=1,2$ and for any bounded interval $[a,b]$.
\end{definition}
\begin{remark}
For smooth solution to (\ref{eq:2}) with (\ref{eq:1-initialdata}), which satisfies the constraint $a^1_{0}-\partial_{x} a^0_1 = 0$, it was proved in \cite{C} that
\[ \partial_{t} A_0 - \partial_{x} A_1 = 0. \]
 Therefore, the Lorentz gauge condition also holds for the strong solution to (\ref{eq:2}) with (\ref{eq:1-initialdata}).
In the remaining, we only consider the Cauchy problem for  (\ref{eq:2}).
\end{remark}

The main result is presented as follows.
\begin{theorem}\label{thm-main}
For the initial data $\psi=(\psi_{1},\psi_{2}) \in L^{2}(R^1)$ and $a_{\mu}=(a^{0}_{\mu},
a^{1}_{\mu}) \in C(R^1)$, there exists a unique global strong solution $(\Psi, A)$ to
(\ref{eq:2}) with (\ref{eq:1-initialdata}), which satisfy
\[
\Psi=(\Psi_{1},\Psi_{2}) \in C([0,\infty);L^{2}(R^1)), A=(A_0,A_1) \in C(R^1\times [0,\infty)).
\]
\end{theorem}

The remaining of the paper is organized as follows. In section 2, we rewrite (\ref{eq:1})
in the equivalent form as (\ref{eq:3}-\ref{eq:5-1}). In section 3, based on the Chadam's result
on the global $H^1$ strong solution for (\ref{eq:1}), we establish the key estimates in
Lemmas \ref{le-conservation} and \ref{le:1} for classical solutions to (\ref{eq:3}-\ref{eq:5-1}), and get the uniform boundness on the solutions. In section 4, the key estimates given in Lemmas \ref{le-conservation} and \ref{le:1} play an important role in proving the precompactness of the approximate solutions $\{(u^n, v^n, A_{\pm}^n)\}_{n=1}^{\infty}$. Then we can get the convergence of the approximate solutions, this yields the existence of the global strong solution there. In section $5$, we prove the uniqueness of the solutions.

\section{Reduction of the problem}
Denote
\begin{eqnarray*}
u=\Psi_1+\Psi_2, v=\Psi_1-\Psi_2, \\
A_+=A_0+A_1, A_-=A_0-A_1,
\end{eqnarray*}
and denote \[u_0=\psi_1+\psi_2, v_0=\psi_1-\psi_2\] and \[a_{\pm}^0=a_0^0\pm a_1^0, a_{\pm}^1=a_0^1\pm a_1^1.\]
Then equations (\ref{eq:2}) can be written equivalently as
\begin{eqnarray}
\partial_{t}u+\partial_{x}u = imv+iA_{+}u, \label{eq:3} \\
\partial_{t}v-\partial_{x}v = imu+iA_{-}v, \label{eq:4}
\\ \Box A_{+} = |v|^2, \label{eq:5}
\\ \Box A_{-} = |u|^2. \label{eq:5-1}
\end{eqnarray}
In addition, the initial data (\ref{eq:1-initialdata}) can be written as
\begin{eqnarray}
u(x,0)=u_0(x), \, v(x,0)=v_0(x) \label{eq:5-inital-uv} \\
A_{\pm}(x,0)=a_{\pm}^0(x), \, \frac{\partial A_{\pm}}{\partial t}(x,0)=a_{\pm}^1(x) \label{eq:5-initial-a}.
\end{eqnarray}

We will get the global solution to (\ref{eq:2}) with (\ref{eq:1-initialdata})
by constructing the solutions to (\ref{eq:3}-\ref{eq:5-initial-a}). That is, we will prove the following,
\begin{theorem}\label{thm-main1}
For any $(u_0,v_0)\in L^2(R^1)$ and $(a_{\pm}^0, a_{\pm}^1)\in C (R^1)$, there exists a unique strong solution $(u,v,A_{\pm})$ to (\ref{eq:3}-\ref{eq:5-initial-a}) with $ (u,v) \in C([0,\infty);L^{2}(R^1))$ and $ A_{\pm} \in C(R^1\times [0,\infty)).$
\end{theorem}
 The strong solution to (\ref{eq:3}-\ref{eq:5-1}) is defined as in Definition \ref{def-strongsolu}.

\section{Some Estimates on the smooth solutions}\label{section-smoothsolu}

In this section, we assume that $\psi_1,\psi_2 \in C_c^{\infty}(R^1)$,
 $a^0_{\nu}\in C^{\infty}_c(R^1)$ and $a^1_{\nu}\in C^{\infty}_c(R^1)$.
 Due to the Chadam's result \cite{C}, the system (\ref{eq:2}) with initial
 data (\ref{eq:1-initialdata}) has a unique solution $(\Psi, A)\in C^{\infty}$
 with $supp \Psi (t,\cdot)$ and $supp A$ are compact sets in $R$ for each $t\in R$.
 This implies that the problem (\ref{eq:3}-\ref{eq:5-initial-a})
 has a unique smooth solution $(u,v,A_{\pm})$.

  We will establish some estimates on the smooth solution $(u,v,A_{\pm})$ in the next. To this end, we
multiply (\ref{eq:3}) and (\ref{eq:4}) by $\overline{u}$ and $\overline{v}$
respectively, which yields the following
\begin{eqnarray}
(|u|^2)_t + (|u|^2)_x = im(\overline{u}v-u\overline{v}),\label{eq:6} \\
(|v|^2)_t - (|v|^2)_x = im(\overline{v}u-v\overline{u}).\label{eq:7}
\end{eqnarray}
Then, direct computation shows that
\begin{lemma}\label{le-conservation}
Suppose $\int^{\infty}_{-\infty} (|u_0(x)|^2+|v_0(x)|^2)dx\le C_0$. Then
\begin{equation}
\int^{\infty}_{-\infty} (|u(x,t)|^2+|v(x,t)|^2)dx=\int^{\infty}_{-\infty} (|u_0(x)|^2+|v_0(x)|^2)dx\le C_0
\end{equation}
for $t\ge 0$.
\end{lemma}

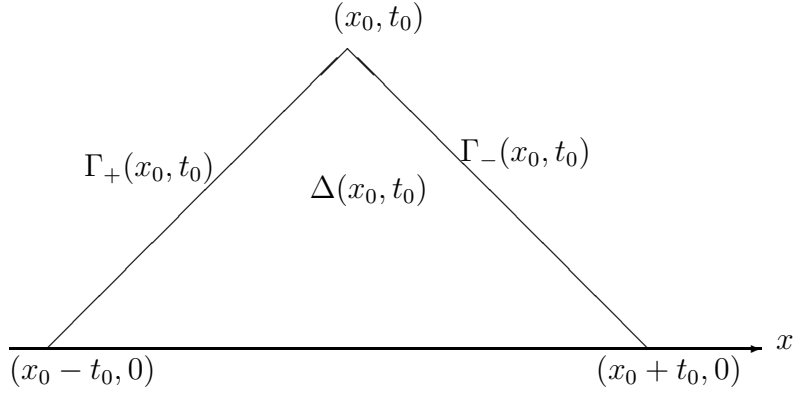
\begin{figure}[h]
\begin{center}
\unitlength=1mm
\begin{picture}(130,75)
\put(40,40){$\Delta(x_0,t_0)$}
\put(10,43){$\Gamma_+(x_0,t_0)$}
\put(60,45){$\Gamma_-(x_0,t_0)$}
\put(85,20){\line(-1,1){40}}
\put(5,20){\line(1,1){40}} \put(43,63){$(x_0,t_0)$}
\put(0,20){\vector(1,0){100}}\put(102,20){$x$}
\put(0,16){$(x_0-t_0, 0)$ } \put(78,16){$(x_0+t_0,0)$}
\end{picture}
\caption{Domain $\Delta(x_0,t_0)$}\label{fig-domain}
\end{center}
\end{figure}

Now, for any $(x_0,t_0)\in R^1\times \overline{R_+}$,
denote \[ \Delta(x_0,t_0) = \{(x,t)|t \in (0,t_0), x \in (x_0 -t_0 +t,x_0 +t_0 -t)\}\]
and
\[ \Gamma_+(x_0,t_0) =\{ (x,t)|t \in [0,t_0], x =x_0 -t_0 +t\},\]
\[ \Gamma_-(x_0,t_0) =\{ (x,t)|t \in [0,t_0], x =x_0 +t_0 -t\},\]
see Fig. \ref{fig-domain}.

\begin{lemma}\label{le:conservation2}
For any $(x_0,t_0)$ with $t_0>0$, there holds the following
\[ \int^{x_0+t_0-t}_{x_0-t_0+t} (|u(x,t)|^2+|v(x,t)|^2)dx \le \int^{x_0+t_0}_{x_0-t_0} (|u_0(x)|^2 +|v_0(x)|^2)dx \] for any $t\in [0,t_0]$.
\end{lemma}
{\it Proof.} For $t\in (0,t_0)$, by (\ref{eq:6}-\ref{eq:7}) we have
\begin{eqnarray*}
(|u|^2+|v|^2)_{t} + (|u|^2-|v|^2)_{x}=0,
\end{eqnarray*}
then,
\begin{eqnarray*}
\frac{d}{dt}\int^{x_0+t_0-t}_{x_0-t_0+t} (|u(x,t)|^2+|v(x,t))|^2)dx &=& \int^{x_0+t_0-t}_{x_0-t_0+t}
\frac{\partial}{\partial t}(|u(x,t)|^2+|v(x,t)|^2)dx \\ &-& (|u(x_0+t_0-t,t)|^2+|v(x_0+t_0-t,t)|^2) \\
&-& (|u(x_0-t_0+t,t)|^2+|v(x_0-t_0+t,t)|^2) \\
&=&-2|u(x_0+t_0-t,t)|^2-2|v(x_0-t_0+t,t)|^2 \\
&\le& 0,
\end{eqnarray*}
which gives the desired result.
The proof is complete.
$\Box$
\begin{lemma}\label{le:1}
  There hold the following
\begin{eqnarray*}
|u(x_0,t_0)|^2  \le q(t_0) \int_{x_0-t_0}^{x_0+t_0} (|u_0(x)|^2+|v_0(x)|^2)dx+ e^{mt_0} |u_0(x_0-t_0)|^2, \\
|v(x_0,t_0)|^2  \le q(t_0) \int_{x_0-t_0}^{x_0+t_0} (|u_0(x)|^2+|v_0(x)|^2)dx+ e^{mt_0} |v_0(x_0+t_0)|^2,
\end{eqnarray*}
where $q(t_0)= e^{mt_0}  m(mt_0+1)$.
\end{lemma}
{\it Proof.}
Integrate the equation (\ref{eq:6}) along its characteristic curve $\Gamma_+(x_0,t_0)$, we have
\[
\frac{d} {d t}|u(x_0-t_0+t,t)|^2  = i m (\overline{u} v - \overline{v} u)(x_0-t_0+t,t).
\]
Then
\[
\frac{d} {d t}|u(x_0-t_0+t,t)|^2 \leq  m(|u(x_0-t_0+t,t)|^2 + |v(x_0-t_0+t,t)|^2),
\]
which leads to
\begin{equation}\label{eq-u-charac1}
|u(x_0-t_0+t,t)|^2 \le {e^{mt_0}} (|u_0(x_0-t_0)|^2 + m \int_{{\Gamma}_+ }|v(x,t)|^2 ds).
\end{equation}
To estimate the righthand side in the above, we
integrate the equation (\ref{eq:7}) over the region $\Delta(x_0,t_0)$, then
\begin{equation}\label{eq-integral-2}
|\int\int_{\Delta(x_0,t_0)}(|v|^2_t - |v|^2_x) dxdt| = |i m \int_0^{t_0} \int_{x_0-t_0+s}^{x_0+t_0-t}(\overline{v}u-\overline{u} v)(y,s) dyds|.
\end{equation}
Applying the divergence theorem and Lemma \ref{le:conservation2} to the lefthand side and righthand side in (\ref{eq-integral-2}) respectively, we have
\[
|2\int_{\Gamma_+(x_0,t_0)} |v|^2 ds+ \int_{x_0-t_0}^{x_0+t_0} |v_0|^2 dx | \leq m t_0 (\int_{x_0-t_0}^{x_0+t_0} (|u_0|^2+|v_0|^2) dx),
\]
which leads to
\begin{equation}\label{eq-v-charac1}
\int_{\Gamma_+(x_0,t_0)}|v|^2 ds \leq (m t_0 +1)(\int_{x_0-t_0}^{x_0+t_0}(|u_0|^2+|v_0|^2) dx).
\end{equation}
With (\ref{eq-u-charac1}) and (\ref{eq-v-charac1}), we have
\[
|u(x_0-t_0+t,t)|^2 \le {e^{mt_0}} m (mt_0+1)\int_{x_0-t_0}^{x_0+t_0} (|u_0|^2+|v_0|^2) dx + {e^{mt_0}} |u_0(x_0-t_0)|^2.
\]
This gives the desired estimate for $|u|^2$. In the same way, we can get the estimate for $|v|^2$.
The proof is complete.
$\Box$

\begin{lemma}\label{le:aymptotic}
For any constant $M\ge 0$,
\begin{eqnarray*}\int_{|y|\ge M}|u(y,\tau)|^2 dy &\le&  2\tau q(\tau) \int_{|x|\ge M-\tau} (|u_0|^2+|v_0|^2)dx \\ &\,& +e^{m\tau} \int_{|y|\ge M} |u_0(y-\tau)|^2 dy,
\end{eqnarray*}
\begin{eqnarray*}
\int_{|y|\ge M}|v(y,\tau)|^2 dy &\le& 2\tau q(\tau)\int_{|x|\ge M-\tau}(|u_0|^2+|v_0|^2)dx
\\ &\,& +e^{m\tau} \int_{|y|\ge M}|v_0(y+\tau)|^2 dy.
\end{eqnarray*}
Here $q(\tau)=e^{m\tau} m(m\tau+1).$
\end{lemma}
{\it Proof.} By Lemma \ref{le:1}, we have,
\begin{eqnarray*}
\int_{|y|\ge M}|u(y,\tau)|^2 dy &\le& \int_{|y|\ge M} q(\tau) \int_{y-\tau}^{y+\tau} (|u_{0}(x)|^2 + |v_{0}(x)|^2) dx dy \\ &\,& + \int_{|y|\ge M} e^{m \tau} |u_{0}(y-\tau)|^2 dy
\\ &\le& e^{m \tau} \int_{|y|\ge M} |u_{0}(y-\tau)|^2 dy
\\ &\,& + q(\tau) \int_{|y|\ge M} \int_{-\tau}^{\tau}(|u_{0}(y+s)|^2 + |v_{0}(y+s)|^2) ds dy \\
&\le& e^{m \tau} \int_{|y|\ge M} |u_{0}(y-\tau)|^2 dy
\\ &\,& + q(\tau) \int_{-\tau}^{\tau} \int_{|x-s|\ge M} (|u_{0}(x)|^2 + |v_{0}(x)|^2) dxds \\
&\le& 2 \tau q(\tau) \int_{|x| \ge M-\tau} (|u_{0}(x)|^2 + |v_{0}(x)|^2) dx \\ &\,& +  e^{m \tau} \int_{|y|\ge M} |u_{0}(y-\tau)|^2 dy.
\end{eqnarray*}
In the same way we can get the estimates on $v$. The proof is complete.
$\Box$

For any $(x_0,t_0)$ and $(x_1,t_1)$ with $t_0\ge 0$ and $t_1\ge 0$, denote
\[ \Omega(x_0,t_0; x_1,t_1)=\Big( \Delta(x_0,t_0)/ \Delta(x_1,t_1)\Big)\cup\Big(\Delta(x_1,t_1)/ \Delta(x_0,t_0)\Big), \] and denote its Lesbeque measure by $meas(\Omega(x_0,t_0; x_1,t_1))$. It is obvious that
\[ \lim_{(x,t)\to (x_0,t_0)} meas(\Omega(x_0,t_0; x,t))=0.\]
\begin{lemma}\label{le:compact-uv}
For any $(x_0,t_0)$ and $(x_1,t_1)$ with $t_0,t_1\in [0,T]$, there hold
\begin{eqnarray*}
 &|{\int\int}_{\Delta(x_0,t_0)}|u(y,\tau)|^2dyd\tau-{\int\int}_{\Delta(x_1,t_1)}|u(y,\tau)|^2dyd\tau |  \\
& \le C_0e^{mT}m(mT+1) meas (\Omega(x_0,t_0; x_1,t_1)) +e^{mT}{\int\int}_{\Omega(x_0,t_0; x_1,t_1)}|u_0(y-\tau)|^2dyd\tau,
\end{eqnarray*}
and
\begin{eqnarray*}
& |{\int\int}_{\Delta(x_0,t_0)}|v(y,\tau)|^2dyd\tau-{\int\int}_{\Delta(x_1,t_1)} |v(y,\tau)|^2dyd\tau |  \\
& \le C_0e^{mT}m(mT+1) meas (\Omega(x_0,t_0; x_1,t_1)) +e^{mT}{\int\int}_{\Omega(x_0,t_0; x_1,t_1)}|v_0(y+\tau)|^2dyd\tau.
\end{eqnarray*}
\end{lemma}
{\it Proof.} First we have
\[
|{\int\int}_{\Delta(x_0,t_0)} |u(y,\tau)|^2dyd\tau-{\int\int}_{\Delta(x_1,t_1)} |u(y,\tau)|^2dyd\tau | = |{\int\int}_{\Omega(x_0,t_0; x_1,t_1)} |u(y,\tau)|^2dyd\tau|.
\]
For the righthand side in above, it follows by Lemma \ref{le:1} that
\begin{eqnarray*}
{\int\int}_{\Omega(x_0,t_0; x_1,t_1)} |u(y,\tau)|^2dyd\tau &\le& {\int\int}_{\Omega(x_0,t_0; x_1,t_1)} e^{m \tau}|u_{0}(y-\tau)|^2dyd\tau \\
&+& {\int\int}_{\Omega(x_0,t_0; x_1,t_1)} q(\tau) \int_{y-\tau}^{y+\tau} (|u_0|^2 + |v_0|^2)dxdyd\tau \\
&\le& e^{mT}{\int\int}_{\Omega(x_0,t_0;x_1,t_1)}|u_{0}(y-\tau)|^2dyd\tau\\ &+& C_{0} q(T) meas(\Omega(x_0,t_0; x_1,t_1))
 ,
\end{eqnarray*}
which proves the first inequality.
In the same way we can prove the second inequality. The proof is complete.
$\Box$

\section{Existence of global strong solution}

To get the strong solution to (\ref{eq:3}-\ref{eq:5-initial-a}) under the assumption of Theorem \ref{thm-main}, we choose sequences of smooth functions $u_0^n, v_0^n, a_{\pm n}^0, a_{\pm n}^1 \in C_c^{\infty}(R^1)$ such that the following hold.
\begin{description}
\item[(A1)] There hold the following
\[\lim\limits_{n\to\infty} \Big(||u_0^n-u_0||_{L^2(R^1)}+||v_0^n-v_0||_{L^2(R^1)}\Big) = 0 \]
and
\[ \lim\limits_{n\to\infty}\Big( ||a_{\pm n}^0-a^0_{\pm}||_{C([a,b])} +||a_{\pm n}^1-a^1_{\pm}||_{C([a,b])} \Big)=0\]
for any bounded interval $[a,b]$.
\item[(A2)]
\[ ||u_0^n||^2_{L^2(R^1)}+||v_0^n||_{L^2(R^1)}^2\le C_0\]
and
\[ \sup\limits_{x\in R} |a_{\pm n}^j(x)|\le 2 C_1\]
for $n\ge 1$ and $j=0,1$, where $C_1=\sup_{x\in R} |a_{+ }^0(x)|+\sup_{x\in R} |a_{+ }^1(x)|+\sup_{x\in R} |a_{- }^0(x)|+\sup_{x\in R} |a_{-}^1(x)|$.
\end{description}

Let $(u^n,v^n,A_+^n, A_-^n)$ be the smooth solutions of (\ref{eq:3}-\ref{eq:5-initial-a}) corresponding to the initial data $(u_0^n, v_0^n, a_{\pm n}^0, a_{\pm n}^1)$. Then, there hold the same estimates  for $(u^n,v^n,A_+^n, A_-^n)$ as in Section \ref{section-smoothsolu}. Moreover,
\begin{equation}\label{eq:wave-1}
 2A_{+}^{n}(x,t) = a_{+ n}^0(x+t) + a_{+ n}^0(x-t) + \int_{x-t}^{x+t} a_{+ n}^1(y) dy +{\int_{0}^{t}} \int_{x-t+s}^{x+t-s} |v^{n}(y,s)|^2 dyds
\end{equation}
and
\begin{equation}\label{eq:wave-2}
 2A_{-}^{n}(x,t) =a_{- n}^0(x+t) + a_{- n}^0(x-t) + \int_{x-t}^{x+t} a_{- n}^1(y) dy +{\int_{0}^{t}} \int_{x-t+s}^{x+t-s} |u^{n}(y,s)|^2 dyds.
\end{equation}

\begin{lemma}\label{le:bdd-a}
For any $T\ge 0$, $n\ge1$, there holds
\[ \sup\limits_{n\ge 0} \sup\limits_{x\in R^1,t\in [0,T]} |A_{\pm}^n(x,t)|\le C_1(T+1)+ C_0T.\]
\end{lemma}
{\it Proof.}
By Lemma \ref{le-conservation}, we can deduce the result from (\ref{eq:wave-1}-\ref{eq:wave-2}). The proof is complete.
$\Box$

\begin{lemma}\label{le:2} There are functions $A^*_{\pm}\in C(R^1\times [0,\infty))$ and subsequences of $\{ A_{\pm}^n\}_{n=1}^{\infty}$, still denoted by $\{ A_{\pm}^n\}_{n=1}^{\infty}$, such that
\[ ||A_+^n-A^*_+||_{C(K)}+ ||A_-^n-A^*_-||_{C(K)} \to 0 \]
as $n\to \infty$ for any compact set $K\subset R^1\times [0,\infty)$.
\end{lemma}
{\it Proof.}
We consider the sequence $\{ A_+^n\}_{n=1}^{\infty}$. Due to the Lemma \ref{le:bdd-a}, it suffices to prove that $ \{\int_{0}^{t} \int_{x-t+s}^{x+t-s} |v^{n}(y,s)|^2 dyds \}_{n=1}^{\infty}$ is equicontinuous in any domain $\Delta (X,T)$ with $T\ge 0$.

By Lemma \ref{le:compact-uv},  for any $(x_0, t_0)$ and $(x_1,t_1)$ with $t_0, t_1\in [0,T]$, we have
\begin{eqnarray}
& |{\int\int}_{\Delta(x_0,t_0)} |v^n(y,\tau)|^2dyd\tau-{\int\int}_{\Delta(x_1,t_1)} |v^n(y,\tau)|^2dyd\tau |  \nonumber\\ & \le C(T) meas (\Omega(x_0,t_0; x_1,t_1)) +C(T){\int\int}_{\Omega(x_0,t_0; x_1,t_1)}|v^n_0(y+\tau)|^2dyd\tau, \nonumber
\\ & \le C(T) meas (\Omega(x_0,t_0; x_1,t_1)) + C(T) {\int\int}_{\Omega(x_0,t_0; x_1,t_1)}|v_0(y+\tau)|^2dyd\tau  \nonumber
\\ & + TC(T)||v_0^n-v_0||^2_{L^2(R^1)}, \label{eq:equicont-v}
\end{eqnarray}
where $C(T)=e^{mT}(m(mT+1)+1)$. Then, by the continuity of the function ${\int\int}_{\Omega(x_0,t_0; x_1,t_1)}|v_0(y+\tau)|^2dyd\tau $ and the convergence of the sequence $\{ v^n\}_{n=1}^{\infty}$, it follows from (\ref{eq:equicont-v}) that $ \{\int_{0}^{t} \int_{x-t+s}^{x+t-s} |v^{n}(y,s)|^2 dyds \}_{n=1}^{\infty}$ is equicontinuous in the domain $\Delta (X,T)$. Therefore the collection $\{ A_+^n\}_{n=1}^{\infty}$ is equicontinuous in the domain $\Delta (X,T)$.

 In the same way, we can prove that $\{ A_-^n\}_{n=1}^{\infty}$ is equicontinuous in the domain $\Delta (X,T)$. Then, by Arzela-Ascoli Theorem (\cite{Rudin}, pp 245), we can have the desired result. The proof is complete.
$\Box$

\begin{lemma}\label{le:convergence-auav}
Let $\{A_{\pm}^{n}\}$ be the sequences given by Lemma $4.2$, For any $T\ge 0$,
\[ \lim\limits_{k,l\to \infty}\int_0^T\int_{R^1} (|A_+^k-A_+^l|^2 + |A_-^k-A_-^l|^2)(|u^k|^2+|v^k|^2)(x,\tau)dxd\tau=0.\]
\end{lemma}
{\it Proof.} For any $M>0$, denote
\[ I_{M,1}^{k,l}=\int_0^T\int_{|x|\ge M} (|A_+^k-A_+^l|^2 + |A_-^k-A_-^l|^2)(|u^k|^2+|v^k|^2)(x,\tau)dxd\tau \]
and
\[ I_{M,2}^{k,l}=\int_0^T\int_{|x|\leq M} (|A_+^k-A_+^l|^2 + |A_-^k-A_-^l|^2)(|u^k|^2+|v^k|^2)(x,\tau)dxd\tau. \]
By Lemmas \ref{le:bdd-a}, \ref{le-conservation} and \ref{le:compact-uv}, we have
\begin{eqnarray*}
  I_{M,1}^{k,l}
&\le & (C_1(T+1)+C_0T) \int_0^T\int_{|x|\ge M} (|u^k|^2+|v^k|^2)(x,\tau)dxd\tau \\
& \le & C^{\prime}_T\Big(T\int_{|x|\ge M-T} (|u_0^k|^2+|v_0^k|^2)dx   +\int_0^T \int_{|y|\ge M} |u_0^k(y-\tau)|^2 dyd\tau \\ &\,& + \int_0^T \int_{|y|\ge M} |v_0^k(y+\tau)|^2 dyd\tau\Big) \\
& \le & C^{\prime}_T\Big(2T\int_{|x|\ge M-T} (|u_0^k|^2+|v_0^k|^2)dx  \Big) \\
&\le& 2TC^{\prime}_T\Big(\int_{|x|\ge M-T} (|u_0|^2+|v_0|^2)dx+||u_0-u_0^k||^2_{L^2}+||v_0-v_0^k||_{L^2}^2  \Big)
\end{eqnarray*}
and
\begin{eqnarray*}
I_{M,2}^{k,l} \le C_0(||A_+^k-A_+^l||_{C(K)}^2 + ||A_-^k-A_-^l||^2_{C(K)}),
\end{eqnarray*}
 where $K=[-M,M]\times[0,T]$ and  $C^{\prime}_T=2e^{mT}(m(mT+1)+1)(C_1(T+1)+C_0T)$. Then, by Lemma \ref{le:2}, we have
 \begin{eqnarray*}
&\,& \limsup\limits_{k,l\to \infty}\int_0^T\int_{R} (|A_+^k-A_+^l|^2 + |A_-^k-A_-^l|^2)(|u^k|^2+|v^k|^2)(x,\tau)dxd\tau \\
&\,& \le \limsup\limits_{k,l\to \infty} (I_{M,1}^{k,l}+I_{M,2}^{k,l}) \\
&\,& \le  2TC^{\prime}_T\Big(\int_{|x|\ge M-T} (|u_0|^2+|v_0|^2)dx\Big),
 \end{eqnarray*}
which implies the desired result by letting $M\to \infty$. The proof is complete. $\Box$

{\bf Proof of the existence of solutions in Theorem \ref{thm-main1}.} Let $\{A_{\pm}^{n}\}$ be the sequences given by Lemma $4.2$. Due to Lemma \ref{le:2}, it suffices to prove the convergence of $\{(u^k,v^k)\}_{k=1}^{\infty} \in (R^1\times (0,T))$, for any $T\ge 0$.

With $(u^n,v^n,A_{\underline{+}}^n)$ of $(2.1)-(2.4)$, we deduce that
\begin{eqnarray*}
(|u^k-u^l|^2)_t + (|u^k-u^l|^2)_x =im((\overline{u}^k-\overline{u}^l)(v^k-v^l)-(u^k-u^l)(\overline{v}^k-\overline{v}^l)) + \\ i(A_+^k-A_+^l)u^{k}(\overline{u}^k-\overline{u}^l) - i(A_+^k-A_+^l)\overline{u}^k(u^k-u^l),\\
(|v^k-v^l|^2)_t - (|v^k-v^l|^2)_x =-im((\overline{u}^k-\overline{u}^l)(v^k-v^l)-(u^k-u^l)(\overline{v}^k-\overline{v}^l)) + \\ i(A_-^k-A_-^l)v^{k}(\overline{v}^k-\overline{v}^l) - i(A_-^k-A_-^l)\overline{v}^k(v^k-v^l).\\
\end{eqnarray*}
Then,  there exists a constant $C_2$ independent of $k,l$ such that
\begin{eqnarray*}
\frac{d}{dt} \int_{R} (|u^k-u^l|^2 + |v^k-v^l|^2) dx \leq C_2\int_{R} (|u^k-u^l|^2 + |v^k-v^l|^2) dx \\
+C_2\int (|A_+^k-A_+^l|^2 + |A_-^k-A_-^l|^2)(|u^k|^2+|v^k|^2)dx
\end{eqnarray*}
for $t\ge 0$.
Applying the Gronwall inequality to this inequality gives
\begin{eqnarray*}
\int_{R} (|u^k-u^l|^2 + |v^k-v^l|^2)(x,t) dx   \leq e^{C_2t}\Big(||u_0^k-u_0^l||_{L^2(R^1)}^2 + ||v_0^v-u_0^l||_{L^2(R^1)}^2  \\
+ C_2 e^{C_2t}\int_0^t\int_{R} (|A_+^k-A_+^l|^2 + |A_-^k-A_-^l|^2)(|u^k|^2+|v^k|^2)(x,\tau)dxd\tau\Big)
\end{eqnarray*}
for $t\ge 0$, which implies $\{(u^k, v^k)\}_{k=1}^{\infty}$ is a Cauchy sequence in $C([0,T], L^2(R^1))$ by Lemma \ref{le:convergence-auav} for any $T\ge 0$. Therefore, there exists a function $(u^*,v^*)\in C([0,T], L^2(R^1))$ such that
\[ ||u^k-u^*||_{C([0,T], L^2(R^1))}+||v^k-v^*||_{C([0,T], L^2(R^1))}\to 0\]
as $k\to\infty$ for any $T\ge 0$. Thus this gives the proof of existence of solution. We will prove the uniqueness in the next section.
$\Box$

\section{Uniqueness of global solutions}

In this section, we prove the uniqueness of the strong solution which we constructed above.

Let $(u_1,v_1,A_{\pm})$ and $(u_2,v_2,B_{\pm})$ be two strong solutions with the same initial data, that is, $u_{1}(x,0)=u_{2}(x,0)$, $v_{1}(x,0)=v_{2}(x,0)$, $A_{\pm}(x,0)=B_{\pm}(x,0)$ and $\partial_{t}A_{\pm}(x,0)=\partial_{t}B_{\pm}(x,0)$. Our goal is to prove that $u_{1}=u_2$, $v_1=v_2$, $A_{\pm}=B_{\pm}$. We assume that
 $(u_1^n,v_1^n,A_{\pm}^n)$ and $(u_2^n,v_2^n,B_{\pm}^n)$ are two sequences of smooth solutions of equations (\ref{eq:3}-\ref{eq:5-1}), such that $(u_1^n,v_1^n,A_{\pm}^n) \to (u_1,v_1,A_{\pm})$, $(u_2^n,v_2^n,B_{\pm}^n) \to (u_2,v_2,B_{\pm})$ in the sense of Definition \ref{def-strongsolu}.

Denote $w_1^n = u_1^n-u_2^n$,  $w_2^n = v_1^n-v_2^n$ and $G_{\pm}^n = A_{\pm}^n-B_{\pm}^n$, then
 \begin{eqnarray*}
 \partial_t w_1^n +\partial_x w_1^n=im w_2^n +G_+^n u_1^n +i B_+^n w_1^n,\\
 \partial_t w_2^n -\partial_x w_2^n=im w_1^n +G_-^n v_1^n +i B_-^n w_2^n,\\
 \Box  G_+^n = \overline{v_1^n} w_2^n + \overline{w_2^n} v_2^n,\\
\Box  G_-^n = \overline{u_1^n} w_1^n + \overline{w_1^n} u_2^n,
 \end{eqnarray*}
which gives the following,
\begin{eqnarray}
\partial_t |w_1^n|^2 +\partial_x |w_1^n|^2 = im(w_2^n \overline{w_1^n}-w_1^n \overline{w_2^n}) + 2 Im(G_+^n \overline{u_1^n} w_1^n), \label{eq:8} \\
\partial_t |w_2^n|^2 -\partial_x |w_2^n|^2 = im(w_1^n \overline{w_2^n} - w_2^n \overline{w_1^n}) + 2 Im(G_-^n \overline{v_1^n} w_2^n),\label{eq:9}\\
\Box  G_+^n = \overline{v_1^n} w_2^n + \overline{w_2^n} v_2^n, \label{eq:10} \\
\Box  G_-^n = \overline{u_1^n} w_1^n + \overline{w_1^n} u_2^n. \label{eq:11}
\end{eqnarray}

By (\ref{eq:8}) and (\ref{eq:9}), we have
\[\partial_{t} (|w_1^n|^2+|w_2^n|^2) +\partial_{x}(|w_1^n|^2-|w_2^n|^2) = 2 Im(G_+^n \overline{u_1^n} w_1^n)+2 Im(G_-^n \overline{v_1^n} w_2^n),\]
for $0 \leq t \leq T$.

Integrating the above equation over the interval $[-T+t,T-t]$,  we have
\begin{eqnarray*}
\frac{d}{dt}({\int_{-T+t}^{T-t}}(|w_1^{n}(x,t)|^2+|w_2^{n}(x,t)|^2)dx)&+& 2|w_1^{n}(T-t,t)|^2+ 2|w_2^{n}(-T+t,t)|^2 \\
&\leq& 2 \int_{-T+t}^{T-t}(|G_+^n \overline{u_1^n} w_1^n| + G_-^n \overline{v_1^n} w_2^n|) dx.
\end{eqnarray*}
Denote $I^{n} (t,T) = \int_{-T+t}^{T-t} (|w_{1}^{n} (x,t)|^2 + |w_{2}^{n} (x,t)|^2) dx.$
Then,
\begin{eqnarray}
\frac{d}{dt} I^{n}(t,T) &\le& 2||G_+^n||_{C(D_T)}||u_1^n||_{L^2(D_T)}||w_1^n||_{L^2(D_T)} {} \nonumber \\
{} &+&  2||G_-^n||_{C(D_T)} ||v_1^n||_{L^2(D_T)} ||w_2^n||_{L^2(D_T)},\label{eq:12}
\end{eqnarray}
where $D_T =\{(x,t)|-T+t < x <T-t,0 \leq t \leq T\}.$

Applying the D'Alembert formula to the equations (\ref{eq:10}-\ref{eq:11}), we have,
\begin{eqnarray*}
2 G_+^n =g_{0+}^n(x-t)+g_{0+}^n(x+t) +\int_{x-t}^{x+t}g_{1+}^n(\tau) d\tau \\ + \int_0^{t} \int_{x-t+s}^{x+t-s} (\overline{v_1^n} w_2^n + \overline{w_2^n} v_2^n)(y,s) dy ds,\\
2 G_-^n =g_{0-}^n(x-t)+g_{0-}^n(x+t) +\int_{x-t}^{x+t}g_{1-}^n(\tau) d\tau \\ \int_0^{t} \int_{x-t+s}^{x+t-s} (\overline{u_1^n} w_1^n + \overline{w_1^n} u_2^n)(y,s) dy ds,
\end{eqnarray*}
where $g_{0\pm}^n(x)=G_{\pm}^n(x,0)$ and $g_{1\pm}^n(x)=\frac{\partial G_{\pm}^n}{\partial t}(x,0)$.

Then,
\begin{eqnarray}
||G_+^{n}||_{C(D_T)} \leq ||g_{0+}^n||_{C([-T,T])} +T||g_{1+}^n||_{(C[-T,T])} \nonumber \\+\frac{1}{2} (||v_1^n||_{L^{2}(D_T)} +||v_2^n||_{L^{2}(D_T)}) ||w_2^n||_{L^{2}(D_T)},\label{eq:13}\\
||G_-^{n}||_{C(D_T)} \leq ||g_{0-}^n||_{C([-T,T])} +T||g_{1-}^n||_{C([-T,T])}\nonumber \\+\frac{1}{2} (||u_1^n||_{L^{2}(D_T)} +||u_2^n||_{L^{2}(D_T)}) ||w_1^n||_{L^{2}(D_T)}. \label{eq:14}
\end{eqnarray}
With (\ref{eq:12}),(\ref{eq:13}) and (\ref{eq:14}), we have
\begin{eqnarray*}
\frac{d}{dt}(I^{n} (t,T)) \leq  C_1(I^{n} (t,T))+ C_2g^n(T)
\end{eqnarray*}
for some constants $C_1>0$ and $C_2>0$ depending only on $C_0$,
 where Lemma \ref{le-conservation} is used. Here $g^n(T)=||g_{0+}^n||_{C([-T,T])}^2 +T^2||g_{1+}^n||_{(C[-T,T])}^2+ ||g_{0-}^n||_{C([-T,T])}^2 +T^2||g_{1-}^n||_{C([-T,T])}^2$.

Then applying the Gronwall inequality and lemma \ref{le-conservation} to above, we can get
\[ I^{n}(t,T) \leq (I^{n}(0,T) +\frac{C_2}{C_1}g^n(T))e^{C_1t}, \]
for $t\ge 0$ and $T\ge 0$.
 Due to the Definition \ref{def-strongsolu},
 \[ \lim_{n\to 0}I^n(0,T)=0, \quad \lim_{n\to 0} g^n(T)=0, \] which lead to
\begin{eqnarray*}
\int_{-T+t}^{T-t}(|u_{1}(y,t)-u_{2}(y,t)|^2 + |v_{1}(y,t)-v_{2}(y,t)|^2)dy &=&\lim_{n \to\infty} I^{n}(t,T)\\
&=&0.
\end{eqnarray*}
Then $u_{1}\equiv v_{1}, u_{2}\equiv v_{2}.$

Now by (\ref{eq:13}), (\ref{eq:14}), we have
\begin{eqnarray*}
||A_{+}-B_{+}||_{C(D_{T})}=\lim_{n\to \infty}||G_{+}^{n}||_{C(D_{T})}=0, \\
||A_{-}-B_{-}||_{C(D_{T})}=\lim_{n\to \infty}||G_{-}^{n}||_{C(D_{T})}=0,
\end{eqnarray*}
therefore $A_{\pm}=B_{\pm}.$
The proof is complete.$\Box$
\section*{ Acknowledgement}

This work was partially  supported by NSFC Project 11031001 and 11121101, by the 111 Project
B08018 and by Ministry of Education of China.

\end{document}